\documentclass[10pt]{amsart}
\usepackage{amsfonts}
\usepackage{ifthen}
\usepackage{amsthm}
\usepackage{amsmath}
\usepackage{graphicx}
\usepackage{amscd,amssymb,amsthm}

\newcounter{minutes}
\divide\time by 60
\newcounter{hours}
\multiply\time by 60 \addtocounter{minutes}{-\time}

\newtheorem{theorem}{Theorem}
\newtheorem{lemma}{Lemma}


\title[\'A. Baricz, R. Sz\'asz/Close-to-convexity of some special functions]{Close-to-convexity of some special functions and their derivatives}

\author[]{\'Arp\'ad Baricz$^{\bigstar}$}
\address{Department of Economics, Babe\c{s}-Bolyai University, 400591 Cluj-Napoca, Romania} \email{bariczocsi@yahoo.com}

\author[]{R\'obert Sz\'asz}
\address{Department of Mathematics and Informatics, Sapientia Hungarian University of Transylvania, 540485 T\^argu-Mure\c{s}, Romania}
\email{rszasz@ms.sapientia.ro}

\thanks{$^{\bigstar}$The research of \'A. Baricz was supported by a research grant of the Romanian National Authority for Scientific Research, CNCS-UEFISCDI, project number PN-II-RU-TE-2012-3-0190/2014.}

\begin{document}

\def\thefootnote{}
\footnotetext{ \texttt{File:~\jobname .tex,
          printed: \number\year-0\number\month-0\number\day,
          \thehours.\ifnum\theminutes<10{0}\fi\theminutes}
} \makeatletter\def\thefootnote{\@arabic\c@footnote}\makeatother

\maketitle

\begin{center}
{\em \'A. Baricz dedicates this paper to the occasion of the 60th birthday of his friend Prof. Tibor K. Pog\'any}
\end{center}

\begin{abstract}
In this paper our aim is to deduce some sufficient (and necessary) conditions for the close-to-convexity of some special functions and their derivatives, like Bessel functions, Struve functions, and a particular case of Lommel functions of the first kind, which can be expressed in terms of the hypergeometric function ${}_1F_2$. The key tool in our proofs is a result of Shah and Trimble about transcendental entire functions with univalent derivatives. Moreover, a known result of P\'olya on entire functions, the infinite product representations and some results on zeros of Bessel, Struve and Lommel functions of the first kind are used in order to achieve the main results of the paper.
\end{abstract}

\section{\bf Introduction and Main Results}

Special functions play an important role in pure and applied mathematics. Bessel functions of the first kind are among of the special functions which were studied by many authors from many different points of view. The geometric properties, like univalence, starlikeness, spirallikeness and convexity were studied already in the sixties by Brown \cite{brown, brown2,brown3}, and also by Kreyszig and Todd \cite{todd}. However, many important problems of Bessel functions, like determining the radius of starlikeness, and the radius of convexity, or finding the optimal parameter for which the normalized Bessel function of the first kind will be starlike, convex, or close-to-convex, have not been studied in details or have not been solved completely. Some of these problems have been studied later in the papers \cite{mathematica,publ,lecture,bsk,samy,basz,szasz,szasz2}, however, there are still some open problems in this direction. For example, there is no information about the close-to-convexity or univalence of the derivatives of Bessel functions, or other special functions. In this paper we make a contribution to the subject by showing some sufficient (and necessary) conditions for the close-to-convexity of some special functions and their derivatives, like Bessel functions, Struve functions, and a particular case of Lommel functions of the first kind, which can be expressed in terms of the hypergeometric function ${}_1F_2$. In order to prove our main results we use a result of Shah and Trimble \cite[Theorem 2]{st} about transcendental entire functions with univalent derivatives. We use also a well-known result of P\'olya on entire functions, and the Weierstrass product representations and some results on zeros of Bessel, Struve and Lommel functions of the first kind are used in order to achieve the main results of the paper. The paper is organized as follows. In this section we recall the result of Shah and Trimble together with the definitions of Bessel, Struve and Lommel functions. Moreover, at the end of this section we present the main results of this paper. Section 2 contains the proofs of these results.

The Bessel function of the first kind $J_{\nu},$ the Struve function of the first kind $\mathbf{H}_{\nu},$ and the Lommel function of the first kind $s_{\mu,\nu},$ are particular solutions of the Bessel differential equation \cite[p. 217]{nist}
\begin{equation}\label{bessel}z^2y''(z)+zy'(z)+(z^2-\nu^2)y(z)=0,\end{equation}
Struve differential equation \cite[p. 288]{nist}
$$z^2y''(z)+zy'(z)+(z^2-\nu^2)y(z)=\frac{\left(\frac{z}{2}\right)^{\nu-1}}{\sqrt{\pi}\Gamma\left(\nu+\frac{1}{2}\right)}$$
and the inhomogeneous Bessel differential equation \cite[p. 294]{nist}
$$z^{2}y''(z) + zy'(z) + (z^{2}-\nu^{2})y(z) = z^{\mu+1}.$$

Let $\mathbb{D}=\{z\in\mathbb{C}:|z|<1\}$ denote the open unit disk. In this paper we are mainly interested on the normalized Bessel function of the first kind $f_{\nu}:\mathbb{D}\to\mathbb{C},$ normalized Struve function of the first kind $h_{\nu}:\mathbb{D}\to\mathbb{C},$ and normalized Lommel functions of the first kind $l_{\mu}:\mathbb{D}\to\mathbb{C},$ which are defined as follows
\begin{equation}\label{fnu}f_{\nu}(z)=2^{\nu}\Gamma(\nu+1)z^{1-\frac{\nu}{2}}J_{\nu}(\sqrt{z})=\sum_{n\geq0}\frac{(-1)^n\Gamma(\nu+1)z^{n+1}}{4^nn!\Gamma(\nu+n+1)},\end{equation}
$$h_{\nu}(z)=\sqrt{\pi}2^{\nu}\Gamma\left(\nu+\frac{3}{2}\right)z^{\frac{1-\nu}{2}}\mathbf{H}_{\nu}(\sqrt{z})=\frac{\sqrt{\pi}}{2}\sum_{n\geq0}\frac{(-1)^n\Gamma\left(\nu+\frac{3}{2}\right)z^{n+1}}{4^n\Gamma\left(n+\frac{3}{2}\right)\Gamma\left(\nu+n+\frac{3}{2}\right)},$$
$$l_{\mu}(z)=\mu(\mu+1)z^{-\frac{\mu}{2}+\frac{3}{4}}s_{\mu-\frac{1}{2},\frac{1}{2}}(\sqrt{z})=
z\cdot{}_1F_2\left(1;\frac{\mu+2}{2},\frac{\mu+3}{2};-\frac{z}{4}\right)=\sum_{n\geq0}\frac{(-1)^n\Gamma\left(\frac{\mu}{2}+1\right)\Gamma\left(\frac{\mu}{2}+\frac{3}{2}\right)z^{n+1}}{4^n\Gamma\left(\frac{\mu}{2}+n+1\right)\Gamma\left(\frac{\mu}{2}+n+\frac{3}{2}\right)},$$
where $J_{\nu}$ and $\mathbf{H}_{\nu}$ stand for the Bessel and Struve functions of the first kind, while $s_{\mu,\nu}$ is the Lommel function of the first kind, which can be expressed in terms of a hypergeometric series as
$$s_{\mu,\nu}(z)=\frac{z^{\mu+1}}{(\mu-\nu+1)(\mu+\nu+1)}{}_{1}F_{2}\left(1;\frac{\mu-\nu+3}{2},\frac{\mu+\nu+3}{2};-\frac{z^2}{4}\right).$$

The next result of Shah and Trimble \cite[Theorem 2]{st} is the cornerstone of this paper.

\begin{lemma}\label{lem1}
Let $f:\mathbb{D}\to\mathbb{C}$ be a transcendental entire function of the form
$$f(z)=z\prod_{n\geq 1}\left(1-\frac{z}{z_n}\right),$$
where all $z_n$ have the same argument and satisfy $|z_n|>1.$ If $f$ is univalent in $\mathbb{D},$ then
\begin{equation}\label{ineq1}
\sum_{n\geq1}\frac{1}{|z_n|-1}\leq 1.
\end{equation}
In fact \eqref{ineq1} holds if and only if $f$ is starlike in $\mathbb{D}$ and all of its derivatives are
close-to-convex there. Furthermore, if $z_n'$ are the zeros of $f'$, then $f$ and all its
derivatives are univalent in $\mathbb{D}$ and map $\mathbb{D}$ onto convex domains if and only if
\begin{equation}\label{ineq2}
\sum_{n\geq1}\frac{1}{|z_n'|-1}\leq 1.
\end{equation}
\end{lemma}

By using the above lemma our aim is to present the following interesting results.

\begin{theorem}\label{th1}
The following assertions are true:
\begin{enumerate}
\item[\bf a.] The function $f_{\nu}$ is starlike and all of its derivatives are close-to-convex in $\mathbb{D}$ if and only if $\nu\geq\nu_0,$ where $\nu_0\simeq-0.5623\dots$ is the unique root of the equation $f'_{\nu}(1)=0$ on $(-1,\infty).$
\item[\bf b.] The function $f_{\nu}$ and all of its derivatives are convex in $\mathbb{D}$ if and only if $\nu\geq\nu_1,$ where $\nu_1\simeq-0.1438\dots$ is the unique root of the equation $3J_{\nu}(1)+2(\nu-2)J_{\nu+1}(1)=0$ on $(-1,\infty).$
\end{enumerate}
\end{theorem}

\begin{theorem}\label{th2}
If $|\nu|\leq \frac{1}{2},$ then $h_{\nu}$ is starlike and all of its derivatives are close-to-convex in $\mathbb{D}.$
\end{theorem}

\begin{theorem}\label{th3}
If $\mu\in\left(-1,1\right),$ $\mu\neq0,$ then $l_{\mu}$ is starlike and all of its derivatives are close-to-convex in $\mathbb{D}.$
\end{theorem}

We note that the results of Theorem \ref{th1} are sharp. Moreover, it is worth to mention that Sz\'asz \cite[Theorem 6]{szasz} deduced already the starlikeness of $f_{\nu},$ while Baricz and Sz\'asz \cite[Theorem 6]{basz} deduced already the convexity of $f_{\nu}$, however, our approach is much easier and as we can see below is applicable also for Struve and Lommel functions. Moreover, in the above theorems we have also information on the close-to-convexity or convexity of the derivatives of the Bessel, Struve and Lommel functions, respectively.

Now, recall that the main idea of this paper is to use Lemma \ref{lem1}, which requires the use of the Weierstrassian infinite canonical representation of the special functions. But, these kind of product representations are not valid for any range of the parameters of the corresponding special functions. Thus, we do not know the best possible range of parameters for which the normalized Struve and Lommel functions of the first kind are starlike and their derivatives will be close-to-convex in the open unit disk. These problems remain open and are subject of further research. As we can see in the case of the Bessel functions of the first kind \cite{bsk,basz,szasz} these kind of problems are not easy to handle, since require a lot of information about the zeros of Bessel functions. But, the zeros of Struve and Lommel functions are not much studied; for example, there is no formula yet for their derivative with respect to the order, which would be an useful source in the study of the geometric properties of Struve and Lommel functions.

\section{\bf Proofs of the Main Results}
\setcounter{equation}{0}

In this section our aim is to present the proofs of the main results.

\begin{proof}[\bf Proof of Theorem \ref{th1}]
{\bf a.} Let $j_{\nu,n}$ denote the $n$th positive zero of the Bessel function $J_{\nu}.$ By using the infinite product representation \cite[p. 235]{nist} of the Bessel functions of the first kind, that is,
$$J_{\nu}(z)=\frac{\left(\frac{z}{2}\right)^{\nu}}{\Gamma(\nu+1)}\prod_{n\geq 1}\left(1-\frac{z^2}{j_{\nu,n}^2}\right),$$
we obtain that the normalized Bessel function of the first kind $f_{\nu}$ belongs to the family of transcendental entire functions, and has the canonical Hadamard factorization
\begin{equation}\label{prod}
f_{\nu}(z)=z\prod_{n\geq1}\left(1-\frac{z}{j_{\nu,n}^2}\right).
\end{equation}
We know \cite[Lemma 3]{szasz} that if $\nu>\nu^{\star},$ where $\nu^{\star}\simeq-0.7745\dots$ is the unique root of $f_{\nu}(1)=0$ or equivalently $j_{\nu,1}=1,$ then
$j_{\nu,1}>1,$ and then we have $j_{\nu,n}>1$ for all $n\in\{1,2,\dots\}$ and $\nu>\nu^{\star}.$ For the sake of completeness we mention here that the above result is almost immediate if use the well-known fact that $\nu\mapsto j_{\nu,1}$ is an increasing function on $(-1,\infty).$ Now, by using \eqref{prod} we obtain that
$$\frac{zf_{\nu}'(z)}{f_{\nu}(z)}=1-\sum_{n\geq1}\frac{1}{j_{\nu,n}^2-z}$$
and hence
\begin{equation}\label{ineq1bes}\frac{f_{\nu}'(1)}{f_{\nu}(1)}=1-\sum_{n\geq1}\frac{1}{j_{\nu,n}^2-1}\geq0\end{equation}
if and only if $\nu\geq\nu_0,$ where $\nu_0\simeq-0.5623\dots$ is the unique root of the equation $f_{\nu}'(1)=0.$ Here we used that when $\nu>\nu^{\star}$ all positive zeros $j_{\nu,n}$ satisfy $j_{\nu,n}>1$ and hence $f_{\nu}(1)>0,$ according to \eqref{prod}. Moreover, we used that for all $\nu>-1$ we have
$$\frac{\partial}{\partial\nu}\left(\frac{f_{\nu}'(1)}{f_{\nu}(1)}\right)=\sum_{n\geq1}\frac{2j_{\nu,n}{\partial j_{\nu,n}}/{\partial\nu}}{(j_{\nu,n}^2-1)^2}\geq0,$$
since the function $\nu\mapsto j_{\nu,n}$ is increasing on $(-1,\infty)$ for all fixed $n\in\{1,2,\dots\},$ see \cite[p. 236]{nist}. Thus, by using the inequality \eqref{ineq1bes} we can see that the function $f_{\nu}$ satisfies \eqref{ineq1} and then by applying Lemma \ref{lem1} we obtain that indeed $f_{\nu}$ is starlike and all of its derivatives are close-to-convex in $\mathbb{D}$ if and only if $\nu\geq\nu_0.$ This completes the proof of this part.

{\bf b.} Differentiating \eqref{fnu} we get
$$f_{\nu}'(z)=2^{\nu-1}\Gamma(\nu+1)z^{-\frac{\nu}{2}}\left((2-\nu)J_{\nu}(\sqrt{z})+\sqrt{z}J_{\nu}'(\sqrt{z})\right).$$
Now, let us denote the $n$th positive zero of the Dini function $z\mapsto (2-\nu)J_{\nu}(z)+zJ_{\nu}'(z)$ by $\beta_{\nu,n}.$ Then the zeros of $f_{\nu}'$ are exactly $\beta_{\nu,n}^2.$ We know that \cite[Lemma 5]{basz} if $\nu>-1,$ then we have
$$\frac{zf_{\nu}''(z)}{f_{\nu}'(z)}=
\frac{\nu(\nu-2)J_\nu(\sqrt{z})+(3-2\nu)\sqrt{z}J_\nu'(\sqrt{z})+
zJ_\nu''(\sqrt{z})}{2(2-\nu)J_\nu(\sqrt{z})+2\sqrt{z}J_\nu'(\sqrt{z})}=-\sum_{n\geq 1}\frac{2z}{\beta_{\nu,n}^2-z}.$$
By using this relation we get
\begin{equation}\label{summa}\sum_{n\geq 1}\frac{1}{\beta_{\nu,n}^2-1}=
-\frac{\nu(\nu-2)J_{\nu}(1)+(3-2\nu)J_{\nu}'(1)+J_{\nu}''(1)}{2(2-\nu)J_{\nu}(1)+2J_{\nu}'(1)}=
\frac{1}{2}\cdot \frac{J_{\nu}(1)+2(1-\nu)J_{\nu+1}(1)}{2J_{\nu}(1)-J_{\nu+1}(1)}.\end{equation}
Here we used the recurrence relation \cite[p. 222]{nist} $$zJ_{\nu}'(z)=\nu J_{\nu}(z)-zJ_{\nu+1}(z),$$ and the fact that the Bessel function $J_{\nu}$ satisfies the Bessel differential equation \eqref{bessel} and thus we have $$z^2J_{\nu}''(z)+zJ_{\nu}'(z)+(z^2-\nu^2)J_{\nu}(z)=0.$$
Now, for $\nu>-1$ let $\gamma_{\nu,n}$ be the $n$th positive root of the equation $\gamma{J}_\nu(z)+zJ_\nu'(z)=0.$ Owing to Landau \cite[p. 196]{landau} we know that if $\nu+\gamma\geq0,$ then the function $\nu\mapsto \gamma_{\nu,n}$ is strictly increasing on $(-1,\infty)$ for $n\in\{1,2,\dots\}$ fixed. This implies that $\nu\mapsto \beta_{\nu,n}$ is strictly increasing on $(-1,\infty)$ for $n\in\{1,2,\dots\}$ fixed, and thus
the function $$\nu\mapsto 1-\sum_{n\geq 1}\frac{1}{\beta_{\nu,n}^2-1}$$ is strictly increasing on $(-1,\infty).$
Consequently we have that $$\sum_{n\geq 1}\frac{1}{\beta_{\nu,n}^2-1}\leq 1$$
if and only if $\nu\geq \nu_1,$ where $\nu_1\simeq-0.1438\dots$ is the unique root of the equation $3J_{\nu}(1)+2(\nu-2)J_{\nu+1}(1)=0$ on $(-1,\infty),$ according to \eqref{summa}. With this the proof is complete.
\end{proof}

\begin{proof}[\bf Proof of Theorem \ref{th2}]
We start with the following result, see \cite[Lemma 1]{bps}. If $|\nu|\leq \frac{1}{2},$ then the Hadamard factorization of the transcendental entire function
$\mathcal{H}_{\nu}:\mathbb{C}\to\mathbb{C},$ defined by
$$\mathcal{H}_{\nu}(z)=\sqrt{\pi}2^{\nu}z^{-\nu-1}\Gamma\left(\nu+\frac{3}{2}\right)\mathbf{H}_{\nu}(z),$$ reads as follows
$$\mathcal{H}_{\nu}(z)=\prod_{n\geq 1}\left(1-\frac{z^2}{h_{\nu,n}^2}\right),$$
where $h_{\nu,n}$ stands for the $n$th positive zero of the Struve function $\mathbf{H}_{\nu}.$ From this we obtain that $h_{\nu}$ belongs to the family of transcendental entire functions, and has the canonical Hadamard factorization
\begin{equation}\label{prod2}
h_{\nu}(z)=z\prod_{n\geq1}\left(1-\frac{z}{h_{\nu,n}^2}\right).
\end{equation}
Recall again that if $\nu>\nu^{\star},$ where $\nu^{\star}\simeq-0.7745\dots$ is the unique root of $f_{\nu}(1)=0$ or equivalently $j_{\nu,1}=1,$ then
$j_{\nu,1}>1,$ and then we have $j_{\nu,n}>1$ for all $n\in\{1,2,\dots\}$ and $\nu>\nu^{\star}.$ Now, since (see Steinig \cite[p. 371]{steinig}) there is exactly one zero of $J_{\nu}$ between two consecutive positive zeros of $\mathbf{H}_{\nu},$ and exactly one in the interval $(0,h_{\nu,1}),$ the above result implies that $h_{\nu,n}>h_{\nu,1}>j_{\nu,1}>1$ for all $|\nu|\leq\frac{1}{2}$ and $n\in\{1,2,\dots\}.$ Thus, we have $h_{\nu}(1)>0$ for all $|\nu|\leq\frac{1}{2}.$ Now, by using \eqref{prod2} we obtain that
$$\frac{zh_{\nu}'(z)}{h_{\nu}(z)}=1-\sum_{n\geq1}\frac{1}{h_{\nu,n}^2-z}$$
and hence
\begin{equation}\label{ineq1Str}\frac{h_{\nu}'(1)}{h_{\nu}(1)}=1-\sum_{n\geq1}\frac{1}{h_{\nu,n}^2-1}\geq0\end{equation}
if and only if $h_{\nu}'(1)\geq0.$ Thus, by using the inequality \eqref{ineq1Str} we can see that the function $h_{\nu}$ satisfies \eqref{ineq1} and then by applying Lemma \ref{lem1} we obtain that indeed $h_{\nu}$ is starlike and all of its derivatives are close-to-convex in $\mathbb{D}$ if and only if $h_{\nu}'(1)\geq0.$

Next, we show that if $\nu\geq -\frac{1}{2},$ then $h_{\nu}'(1)>0.$ For this we use the recurrence relation
$$\mathbf{H}_{\nu-1}(x)=\frac{\nu}{x}\mathbf{H}_{\nu}(x)+\mathbf{H}_{\nu}'(x)$$
and we obtain
$$2zh_{\nu}'(z^2)=\sqrt{\pi}2^{\nu}\Gamma\left(\nu+\frac{3}{2}\right)z^{-\nu}\left((1-2\nu)\mathbf{H}_{\nu}(z)+z\mathbf{H}_{\nu-1}(z)\right),$$
which implies in particular
$$h_{\nu}'(1)=\sqrt{\pi}2^{\nu}\Gamma\left(\nu+\frac{3}{2}\right)\left((1-2\nu)\mathbf{H}_{\nu}(1)+\mathbf{H}_{\nu-1}(1)\right).$$
By using \cite[p. 291]{nist}
$$\mathbf{H}_{-\frac{1}{2}}(z)=\sqrt{\frac{2}{\pi z}}\sin z\ \ \ \ \mbox{and}\ \ \ \ \mathbf{H}_{-\frac{3}{2}}(z)=\sqrt{\frac{2}{\pi z}}\left(\cos z-\frac{\sin z}{z}\right)$$
for $\nu=-\frac{1}{2}$ we have
$$(1-2\nu)\mathbf{H}_{\nu}(1)+\mathbf{H}_{\nu-1}(1)=\sqrt{\frac{2}{\pi}}\left(\sin 1+\cos 1\right)\simeq1.102495575\dots,$$
that is, we have $h_{-\frac{1}{2}}'(1)>0.$ Now, by using the recurrence relation \cite[p. 292]{nist}
$$\mathbf{H}_{\nu-1}(x)+\mathbf{H}_{\nu+1}(x)=\frac{2\nu}{x}\mathbf{H}_{\nu}(x)+
\frac{\left(\frac{x}{2}\right)^{\nu}}{\sqrt{\pi}\Gamma\left(\nu+\frac{3}{2}\right)}$$
and the integral representation of $\mathbf{H}_{\nu}$
$$\mathbf{H}_{\nu}(x)=\frac{2\left(\frac{x}{2}\right)^{\nu}}{\sqrt{\pi}\Gamma\left(\nu+\frac{1}{2}\right)}
\int_0^1(1-t^2)^{\nu-\frac{1}{2}}\sin(xt)dt,$$
we get that when $\nu>-\frac{1}{2}$
\begin{align*}h_{\nu}'(1)&=\sqrt{\pi}2^{\nu}\Gamma\left(\nu+\frac{3}{2}\right)\left(\mathbf{H}_{\nu}(1)-\mathbf{H}_{\nu+1}(1)+
\frac{1}{\sqrt{\pi}2^{\nu}\Gamma\left(\nu+\frac{3}{2}\right)}\right)\\
&=1+\int_0^1(1-t^2)^{\nu-\frac{1}{2}}(2\nu+t^2)\sin(t)dt>1-\int_0^1(1-t^2)^{\nu+\frac{1}{2}}\sin(t)dt>0,\end{align*}
since the last integrand is less than $1.$ This completes the proof.
\end{proof}

\begin{proof}[\bf Proof of Theorem \ref{th3}]
We start with the following result, see \cite[Lemma 1]{koumandos}. Let
$$\varphi_{k}(z)={}_{1}F_{2}\left(1;\,\frac{\mu-k+2}{2},\frac{\mu-k+3}{2};-\frac{z^2}{4}\right),$$
where $z\in\mathbb{C},$ $\mu\in\mathbb{R}$ and $k\in\{0,1,\ldots\}$ such that $\mu-k$ is not in $\{0,-1,\dots\}.$ Then, the Hadamard's factorization of $\varphi_{k}$ is of the form
\begin{equation}\label{prodlo}
\varphi_{k}(z)=\prod_{n\geq1}\left(1-\frac{z^2}{z_{\mu,k,n}^2}\right),\end{equation}
where $\pm z_{\mu,k,1}, \pm z_{\mu,k,2},\ldots$ are all zeros of the function $\varphi_{k}$ and the infinite product is absolutely convergent.
For $n\in\{1,2,\dots\}$ let $\xi_{\mu,n}:=z_{\mu,0,n}$ be the $n$th positive zero of $\varphi_{0},$ and let $\xi_{\mu,0}=0$. Then, by using \eqref{prodlo} for $k=0$ we get that $$l_{\mu}(z)=z\varphi_0(\sqrt{z})=z\prod_{n\geq1}\left(1-\frac{z}{\xi_{\mu,n}^2}\right)$$ satisfies
$$\frac{l_{\mu}'(z)}{l_{\mu}(z)}=\frac{1}{z}-\sum_{n\geq 1}\frac{1}{\xi_{\mu,n}^2-z},$$
which implies that
\begin{equation}\label{lomm}\frac{l_{\mu}'(1)}{l_{\mu}(1)}=1-\sum_{n\geq 1}\frac{1}{\xi_{\mu,n}^2-1}.\end{equation}
We know that (see \cite[Lemma 2.1]{lamp}) $\xi_{\mu,n}\in\left(n\pi,(n+1)\pi\right)$ for all $\mu\in(0,1)$ and $n\in\{1,2,\dots\},$ which implies that $\xi_{\mu,n}>\xi_{\mu,1}>\pi>1$ for all $\mu\in(0,1)$ and $n\in\{1,2,\dots\}.$ Consequently, we have that $l_{\mu}(1)=\varphi_0(1)>0$ for $\mu\in(0,1).$ On the other hand, differentiating both sides of $l_{\mu}(z^2)=z^2\varphi_0(z)$ we obtain $2zl_{\mu}'(z^2)=2z\varphi_0(z)+z^2\varphi_0'(z).$ Now, by using the recurrence relation \cite[Lemma 2]{koumandos} $$(\mu+1)\varphi_1(z)=(\mu+1)\varphi_0(z)+z\varphi_0'(z)$$
and the integral representations \cite[Lemma 3]{koumandos}
$$z\varphi_{0}(z)=\mu(\mu+1)\int_{0}^{1}(1-t)^{\mu-1}\sin(zt)dt,\ \ \varphi_{1}(z)=\mu\int_{0}^{1}(1-t)^{\mu-1}\cos(zt)dt,$$
we get for $\mu\in(0,1)$ that
$$2l_{\mu}'(1)=(\mu+1)\varphi_1(1)+(1-\mu)\varphi_0(1)=\mu(\mu+1)\int_0^1(1-t)^{\mu-1}\left(\cos(t)+(1-\mu)\sin(t)\right)dt>0.$$
Thus, for $\mu\in(0,1)$ the right-hand side of \eqref{lomm} is positive, and using Lemma \ref{lem1} we get that for $\mu\in(0,1)$ the transcendental entire function $l_{\mu}$ is starlike in $\mathbb{D}$ and all of its derivatives are close-to-convex in $\mathbb{D}.$

Now, let $\zeta_{\mu,n}:=z_{\mu,1,n}$ be the $n$th positive zero of $\varphi_{1}$. By applying \eqref{prodlo} for $k=1$ we get that
$$l_{\mu-1}(z)=z\varphi_1(\sqrt{z})=z\prod_{n\geq1}\left(1-\frac{z}{\zeta_{\mu,n}^2}\right)$$ satisfies
$$\frac{l_{\mu-1}'(z)}{l_{\mu-1}(z)}=\frac{1}{z}-\sum_{n\geq 1}\frac{1}{\zeta_{\mu,n}^2-z},$$
which implies that
\begin{equation}\label{lomm2}\frac{l_{\mu-1}'(1)}{l_{\mu-1}(1)}=1-\sum_{n\geq 1}\frac{1}{\zeta_{\mu,n}^2-1}.\end{equation}
Recall the following known result of P\'olya \cite{pol}: if the function $f$ is positive, strictly increasing, continuous on $[0,1)$ and satisfies $\int_0^1f(t)dt<\infty,$ then the entire function $z\mapsto \int_0^1f(t)\cos(zt)dt$ has only real and simple zeros, and each interval $\left((2n+1)\frac{\pi}{2},(2n+3)\frac{\pi}{2}\right),$ $n\in\{1,2,\dots\},$ contains exactly one zero. Applying this result we obtain that $\zeta_{\mu,n}>\zeta_{\mu,1}>\frac{\pi}{2}>1$ for all $n\in\{1,2,\dots\},$ $\mu\in(0,1),$ which in turn implies that $l_{\mu-1}(1)>0$ for $\mu\in(0,1).$
Now, differentiating both sides of $l_{\mu-1}(z^2)=z^2\varphi_1(z)$ we obtain $2zl_{\mu-1}'(z^2)=2z\varphi_1(z)+z^2\varphi_1'(z),$
which by means of the integral representation of $\varphi_1$ implies that
$$2l_{\mu-1}'(1)=2\varphi_1(1)+\varphi_1'(1)=\mu\int_0^1(1-t)^{\mu-1}\left(2\cos(t)-t\sin(t)\right)dt>0,$$
where $\mu\in(0,1).$ Here we used that the function $g:[0,1]\to \mathbb{R},$ defined by $g(t)=2\cos(t)-t\sin(t),$ is decreasing and hence for all $t\in[0,1]$ we have $g(t)\geq g(1)=2\cos 1-\sin 1\simeq0.2391336269\dots>0.$ Thus, for $\mu\in(0,1)$ the right-hand side of \eqref{lomm2} is positive, and by using Lemma \ref{lem1} we conclude that for $\mu\in(0,1)$ the transcendental entire function $l_{\mu-1}$ is starlike in $\mathbb{D}$ and all of its derivatives are close-to-convex in $\mathbb{D}.$ Changing $\mu$ to $\mu+1,$ the proof of this theorem is complete.
\end{proof}

\end{document}